\documentclass{amsart}

\usepackage{amsmath,amscd,amssymb,amsthm}
\usepackage{latexsym,enumerate,graphicx,geometry,rotating,braket,mathrsfs,url,graphicx}
\usepackage[all]{xy}
\usepackage[sort&compress,numbers]{natbib}
\usepackage[pdftex]{hyperref}
\usepackage{extarrows}
\usepackage{stmaryrd}
\usepackage{subfigure}
\usepackage{setspace}
\usepackage{float}

\numberwithin{equation}{section}

\theoremstyle{plain}
\newtheorem{thm_}[equation]{Theorem}

\newtheorem{lemma_}[equation]{Lemma}
\newtheorem{prop_}[equation]{Proposition}
\newtheorem{cor_}[equation]{Corollary}
\newtheorem{eg_}[equation]{Example}
\newtheorem{con_}[equation]{Conjecture}
\newtheorem*{cons_}{Conjecture}

\theoremstyle{definition}
\newtheorem{thmu_}[equation]{Theorem}
\newtheorem*{thmus_}{Theorem}
\newtheorem{propu_}[equation]{Proposition}
\newtheorem*{propus_}{Proposition}
\newtheorem{coru_}[equation]{Corollary}
\newtheorem{lemu_}[equation]{Lemma}

\newtheorem*{egus_}{Example}
\newtheorem{def_}[equation]{Definition}
\newtheorem*{defs_}{Definition}
\newtheorem{ex_}[equation]{Example}

\newcommand{\thm}[1]{\begin{thm_}#1\end{thm_}}

\newcommand{\lemm}[1]{\begin{lemma_}#1\end{lemma_}}

\newcommand{\prop}[1]{\begin{prop_}#1\end{prop_}}

\newcommand{\defi}[1]{\begin{def_}#1\end{def_}}

\newcommand{\pf}[1]{\begin{proof}#1\end{proof}}

\DeclareMathOperator{\rank}{rank}

\newcommand{\FF}{\mathbb F}

\newcommand{\QQ}{\mathbb Q}

\newcommand{\ZZ}{\mathbb Z}

\renewcommand{\mod}[1]{\ \mathrm{mod}\ #1}  

\usepackage[OT2,T1]{fontenc}

\DeclareSymbolFont{cyrletters}{OT2}{wncyr}{m}{n}

\DeclareMathSymbol{\Sha}{\mathalpha}{cyrletters}{"58}

\begin{document}

\title{Constant Tamagawa numbers of special elliptic curves}
\author[L. Li]{Luying Li}
\address{State Key Laboratory of Information Security\\
Institute of Information Engineering\\
Chinese Academy of Sciences\\
Beijing 100093, P.R. China}
\email{lilouying16@mails.ucas.ac.cn}

\author[C. Lv]{Chang Lv}
\address{State Key Laboratory of Information Security\\
Institute of Information Engineering\\
Chinese Academy of Sciences\\
Beijing 100093, P.R. China}
\email{lvchang@amss.ac.cn}

\thanks{This work was supported by National Natural Science Foundation of China (Grant No.11701552)}

\maketitle


\begin{abstract}
For the elliptic curves $E_{\sigma 2D} : y^2 = x^3 + \sigma 2Dx$ , which has 2-isogeny curve $E'_{\sigma 2D} : y^2 = x^3 -\sigma 8Dx$,  $\sigma = \pm 1,\ D = p_1^{e_1}p_2^{e_2}\cdots p_n^{e_n}$, where $p_i$ are different odd prime numbers and $e_i = 1 \text{ or } 3$, we demonstrate that Tamagawa numbers of these elliptic curves are always one or zero by the use of matrix in finite field $\FF_2$. The specific number depends on the value of $\sigma$. By our proofs of these results, we find a method to quickly sieve a part of the elliptic curves with Mordell-Weil rank zero or rank one in this form as an application.
\keywords{Elliptic curves, Selmer groups, Tamagawa numbers}
\end{abstract}

\section{Introduction}
The Tamagawa numbers for $E_A: y^2 = x^3-A x$ with $A \in \QQ^\times$ are
$$T_{A} = \# S^{(\phi)}(E_A)/\# S^{(\phi)}(E_{-4A}) = 2^{-t_A}.$$ 
$S^{(\phi)}(E_A)$ is the $\phi$-Selmer group of $E_A$, $\phi$ is the morphism from $E_A$ to its 2-isogeny has the following form.
$$
\phi : E \rightarrow E', \quad (x,y)\mapsto (\frac{x^2}{y^2},\frac{y(-A-x^2)}{x^2})
$$

In this paper, we mainly concern the $\phi$-Selmer groups $S^{({\phi})}({E_{\pm}})$ of following elliptic curves:
$$
E = E_{-\sigma2D}: y^2 = x^3 + \sigma2D x, \ \sigma = \pm 1, \  D = p_1^{e_1}p_2^{e_2}\cdots p_n^{e_n} \in \QQ^\times ,\ e_i = 1,3.
$$
We show the Tamagawa numbers $t_{-\sigma 2 D}$ in this subgroup are uniform distribution of $\{0,1\}$.

Let $E' = E'_{\sigma 2D}: y^2 = x^3 -\sigma 8D x$ denote the 2-isogeny of $E$, and let $E_+ = y^2 = x^3 + 2 D x$, $E_+' = y^2 = x^3 - 8 D x$ and $E_- = y^2 = x^3 - 2 D x$, $E_-' = y^2 = x^3 + 8 D x$ for simplification.
Our main purpose is to prove the following theorems:
\thm{[Theorem \ref{thm_n1}] \label{thm_1} Let $E_-$, $E_-'$ be the elliptic curves defined above. Then     $$2^{-t_{2D}} = T_{2D} = \# S^{(\phi)}(E_-)/\# S^{(\phi)}(E'_{-}) = 1/2.$$ 
}
\thm{[Theorem \ref{thm_n2}] \label{thm_2} Let $E_+$, $E_+'$ be the elliptic curves defined above. Then     $$2^{-t_{-2D}} = T_{-2D} = \# S^{(\phi)}(E_+)/\# S^{(\phi)}(E'_{+}) = 1.$$
}

Furthermore, we claim that for all D, the parity of $\dim_2 S_2(E_-)$ is odd and the parity of $\dim_2 S_2(E_+)$ is even.
The definition of $\phi$-Selmer groups can be found in, for example \cite{gtm106}. From the proof of our main theory, we find a method to quickly sieve a part of the elliptic curves with Mordell-Weil rank zero or rank one in this form.

For $D$ which has quite few prime factors, we can enumerate all the possible cases of $S^{(\phi)}(E_\pm)$ and $S^{({\phi})}({E_\pm'})$ to verify the above two theorems. For example, if $D = p_{1}p_{2}$, we have Table 1 for $\sigma = -1$, and Table 2 for $\sigma = 1$, from which, it is easy to obtain our theorems.

\begin{table}[H]
\begin{center}

\begin{tabular}{cccc|l|l}
	$p \mod 8$ & $q \mod 8$ & $(\frac{p}{q})$ &($\frac{q}{p})$ &$S^{(\phi)}(E_-)$ &$S^{({\phi})}({E_-'})$ \\ \hline
	1 & 1 & -1 & -1 & 1,pq,2,2pq& 1,-1,pq,-pq,2,-2,2pq,-2pq\\
	1 & 1 & 1 & 1 & 1,p,q,pq,2,2p,2q,2pq & 1,-1,p,-p,q,-q,pq,-pq,2,-2,2p,-2p,2q,-2q,2pq,-2pq\\
	1 & 3 & -1 & -1 & 1,2pq & 1,pq,-2,-2pq\\
	1 & 3 & 1 & 1 & 1,p,2q,2pq & 1,p,q,pq,-pq,-2,-2p,-2q,-2pq\\
	1 & 5 & -1 & -1 & 1,2pq & 1,-1,2pq,-2pq\\
	1 & 5 & 1 & 1 & 1,p,2q,2pq & 1,-1,p,-p,2q,-2q,2pq,-2pq\\
	1 & 7 & -1 & -1 & 1,2pq & 1,-pq,2,-2pq\\
	1 & 7 & 1 & 1 & 1,p,2q,2pq & 1,p,-q,-pq,2,2p,-2q,-2pq\\
	3 & 1 & -1 & -1 & 1,2pq & 1,pq,-2,-2pq\\
	3 & 1 & 1 & 1 & 1,q,2p,2pq & 1,p,q,pq,-2,-2p,-2q,-2pq\\
	3 & 3 & -1 & 1 & 1,2pq &1,pq,-2,-2pq\\
	3 & 3 & 1 & -1 & 1,2pq &1,pq,-2,-2pq\\
	3 & 5 & -1 & -1 & 1,2pq &1,-q,2p,-2pq\\
	3 & 5 & 1 & 1 & 1,2pq &1,p,-2q,-2pq\\
	3 & 7 & -1 & 1 & 1,2pq &1,q,-2p,-2pq\\
	3 & 7 & 1 & -1 & 1,2pq &1,-q,2p,-2pq\\
	5 & 1 & -1 & -1 & 1,2pq &1,-1,2pq,-2pq\\
	5 & 1 & 1 & 1 & 1,q,2p,2pq &1,-1,q,-q,2p,-2p,2pq,-2pq\\
	5 & 3 & -1 & -1 & 1,2pq &1,-p,2q,-2pq\\
	5 & 3 & 1 & 1 & 1,2pq &1,q,-2p,-2pq\\
	5 & 5 & -1 & -1 & 1,2pq &1,-1,2pq,-2pq \\
	5 & 5 & 1 & 1 & 1,2pq &1,-1,2pq,-2pq\\
	5 & 7 & -1 & -1 & 1,2pq &1,-p,2q,-2pq\\
	5 & 7 & 1 & 1 & 1,2pq &1,-q,2p,-2pq\\
	7 & 1 & -1 & -1 & 1,2pq &1,-pq,2,-2pq\\
	7 & 1 & 1 & 1 & 1,q,2p,2pq &1,q,-p,-pq,2,-2p,2q,-2pq\\
	7 & 3 & -1 & 1 & 1,2pq &1,-p,2q,-2pq\\
	7 & 3 & 1 & -1 & 1,2pq &1,p,-2q,-2pq\\
	7 & 5 & -1 & -1 & 1,2pq &1,-q,2p,-2pq\\
	7 & 5 & 1 & 1 & 1,2pq &1,-p,2q,-2pq\\
	7 & 7 & -1 & 1 & 1,pq,2,2pq &1,q,-p,-pq,2,-2p,2q,-2pq\\
	7 & 7 & 1 & -1 & 1,pq,2,2pq &1,p,-q,-pq,2,2p,-2q-2pq\\
	\hline
\end{tabular}
	\caption{table 1}
\end{center}
\end{table}
\newpage
\begin{table}[h]
\begin{center}
\begin{tabular}{cccc|l|l}
	$p \mod 8$ & $q \mod 8$ & $(\frac{p}{q})$ &($\frac{q}{p})$ &$S^{({\phi})}({E_+'})$ &$S^{(\phi)}(E_+)$ \\ \hline
	1 & 1 & -1 & -1 & 1,pq,2,2pq& 1,pq,-2,-2pq\\
	1 & 1 & 1 & 1 & 1,p,q,pq,2,2p,2q,2pq & 1,p,q,pq,-2,-2p,-2q,-2pq\\
	1 & 3 & -1 & -1 & 1,2pq & 1,-2pq\\
	1 & 3 & 1 & 1 & 1,p,2q,2pq & 1,p,-2q,-2pq\\
	1 & 5 & -1 & -1 & 1,2pq & 1,-2pq\\
	1 & 5 & 1 & 1 & 1,p,2q,2pq & 1,p,-2q,-2pq\\
	1 & 7 & -1 & -1 & 1,pq,2,2pq & 1,-pq,2,-2pq\\
	1 & 7 & 1 & 1 & 1,p,q,pq,2,2p,2q,2pq & 1,p,-q,-pq,2,2p,-2q,-2pq\\
	3 & 1 & -1 & -1 & 1,2pq & 1,-2pq\\
	3 & 1 & 1 & 1 & 1,q,2p,2pq & 1,q,-2p,-2pq\\
	3 & 3 & -1 & 1 & 1,q,2p,2pq &1,pq,-2,-2pq\\
	3 & 3 & 1 & -1 & 1,p.2q,2pq &1,pq,-2,-2pq\\
	3 & 5 & -1 & -1 & 1,2pq &1,-2pq\\
	3 & 5 & 1 & 1 & 1,2pq &1,-2pq\\
	3 & 7 & -1 & 1 & 1,2pq &1,-2pq\\
	3 & 7 & 1 & -1 & 1,p,2q,2pq &1,-q,2p,-2pq\\
	5 & 1 & -1 & -1 & 1,2pq &1,-2pq\\
	5 & 1 & 1 & 1 & 1,q,2p,2pq &1,q,-2p,-2pq\\
	5 & 3 & -1 & -1 & 1,2pq &1,-2pq\\
	5 & 3 & 1 & 1 & 1,2pq &1,-2pq\\
	5 & 5 & -1 & -1 & 1,2pq &1,-2pq \\
	5 & 5 & 1 & 1 & 1,2pq &1,-2pq\\
	5 & 7 & -1 & -1 & 1,2pq &1,-2pq\\
	5 & 7 & 1 & 1 & 1,q,2p,2pq &1,-q,2p,-2pq\\
	7 & 1 & -1 & -1 & 1,pq,2,2pq &1,-pq,2,-2pq\\
	7 & 1 & 1 & 1 & 1,p,q,pq,2,2p,2q,2pq &1,q,-p,-pq,2,-2p,2q,-2pq\\
	7 & 3 & -1 & 1 & 1,q,2p,2pq &1,-p,2q,-2pq\\
	7 & 3 & 1 & -1 & 1,2pq &1,-2pq\\
	7 & 5 & -1 & -1 & 1,2pq &1,-2pq\\
	7 & 5 & 1 & 1 & 1,p,2q,2pq &1,-p,2q,-2pq\\
	7 & 7 & -1 & 1 & 1,pq,2,2pq &1,-p,2q,-2pq\\
	7 & 7 & 1 & -1 & 1,pq,2,2pq &1,-q,2p,-2pq\\
	\hline
\end{tabular}
	\caption{table 2}
\end{center}
\end{table}

The paper is organized as follows: in Section 2, we introduce some basic facts from the literature; in Section 3 we prove Theorems \ref{thm_1} and \ref{thm_2} by using the matrix in $\FF_2$. Finally, in Section 4, we show a method
which can quickly sieve a part of the elliptic curves with Mordell-Weil rank zero or rank one in this form.
\section{Notations and Lemmas}\label{sec_2}

\[
\begin{tabular}{l|l}
\text{Symbols} & \text{Meanings} \\ \hline
$D$  &  $D = p_1^{e_1}p_2^{e_2}\cdots p_n^{e_n}$, $p_i$ are different odd prime numbers and $e_i = 1,3$\\
$n$  &  number of odd prime factors of $D$\\
$p_{n+1}$ & $p_{n+1} = -1$ \\
$S$ & $S = \{ \infty, 2, p_i \}$ \\
$\mathbb{Q}(S,2)$ & $\mathbb{Q}(S,2) = <-1,2,p_i>\subseteq \mathbb{Q}^{\times}/\QQ^{\times2} $, $p_i$ is prime factor of D  \\
$\QQ_{2D}$ & $\QQ_{2D} = <p_i> \subseteq \QQ^{\times}/\QQ^{\times2} $\\
$\QQ_{-2D}$ & $\QQ_{-2D} = <-1, p_i> \subseteq \QQ^{\times}/\QQ^{\times2}$\\
$d$ & an element of $\mathbb{Q}(S,2)$, usually write as $d = p_{i_1}p_{i_2}\cdots p_{i_k} , 1 \le i_1 < i_2 < \cdots < i_k \le n+1$ \\

$\delta (d)$ & the sum of $y_{i_j}(-1)$, $p_{i_j}$ are different prime factors of $d$. The definition of $y_{i_j}(-1)$ can be found at \ref{eq:array2}\\
$(\frac{\ }{\ })$ & the Legendre  Symbol\\
\hline
\end{tabular}
\]
For $d\in \QQ(S,2)$, define the following curves: $ C_{d-}: W^2 = d + \frac{8DZ^4}{d},$ $ C'_{d-}: W^2 = d - \frac{2DZ^4}{d},$ $ C_{d+}: W^2 = d - \frac{8DZ^4}{d},$ $ C'_{d+}: W^2 = d + \frac{2DZ^4}{d}. $
According to \cite{gtm106}, we have
$$
\begin{aligned}
 S^{(\phi)}(E_-) \cong \{d \in \mathbb{Q}(S,2): C_{d-}(\mathbb{Q}_v) \neq \emptyset , \forall v \in S \} \subseteq \QQ^{\times}/\QQ^{\times2 }, \\
S^{({\phi})}(E_-') \cong \{d \in \mathbb{Q}(S,2): C_{d-}'(\mathbb{Q}_v) \neq \emptyset , \forall v \in S \}\subseteq \QQ^{\times}/\QQ^{\times2 }. \\
S^{(\phi)}(E_+) \cong \{d \in \mathbb{Q}(S,2): C_{d+}(\mathbb{Q}_v) \neq \emptyset , \forall v \in S \} \subseteq \QQ^{\times}/\QQ^{\times2 },\\
S^{({\phi})}(E_+') \cong \{d \in \mathbb{Q}(S,2): C_{d+}'(\mathbb{Q}_v) \neq \emptyset , \forall v \in S \}\subseteq \QQ^{\times}/\QQ^{\times2 }, \\
\end{aligned}
$$
and we have the following inclusions:
$\{1, \ 2D\} \subseteq S^{(\phi)}(E_-) $, $\{1, -2D\} \subseteq S^{({\phi})}({E_-'})$ and $\{1, \ -2D\} \subseteq S^{(\phi)}(E_+) $, $\{1, 2D\} \subseteq S^{({\phi})}({E_+'})$. Since $S^{(\phi)}(E_\pm) $ and $S^{({\phi})}({E_\pm'})$ have above properties, and are
multiplicative abelian groups, when we calculate $S^{(\phi)}(E_\pm)$ and $S^{({\phi})}({E_\pm'})$,
we only need to concern ourselves with $d \in \QQ_{2D}$ and $\QQ_{-2D}$.

\defi{Let

\begin{itemize}
\item $s_-(\phi) \triangleq  S^{(\phi)}(E_-) \bigcap \QQ_{2D}$,
\item $s_-'(\phi) \triangleq S^{({\phi})}({E_-'}) \bigcap \QQ_{-2D}$,
\item $s_+(\phi) \triangleq  S^{(\phi)}(E_+) \bigcap \QQ_{-2D}$,
\item $s_+'(\phi) \triangleq S^{({\phi})}({E_+'}) \bigcap \QQ_{2D}$.
\end{itemize}
}
Apparently, there is a nature map
$$
\begin{aligned}
s_-(\phi) &\rightarrow S^{(\phi)}(E_-) - s_-(\phi)\\
d &\mapsto d \times 2D (\mod \QQ^{\times}/\QQ^{\times 2})  \\
\end{aligned}
$$
and it is easy to check that it is a one-to-one map. Thus we have $ \#S^{(\phi)}(E_-) = 2\# s_-(\phi)$. Similarly, $ \# S^{(\phi)}(E_+) = 2\# s_+(\phi)$, $\#S^{({\phi})}({E_\pm'}) = 2 \# s_\pm'(\phi)$.

The following lemmas give the method to calculate $ s_\pm(\phi)$ and $s_\pm'(\phi) $ by using the Hensel's Lemma. The proof of these lemmas can be found in \cite{zeng}.

\lemm{\label{lem1}
Let $C_{d-}$ be defined as above. Then\vspace{1.5ex}

\begin{enumerate}
\item[(1)] $C_{d-}(\mathbb{Q}_2) \neq \emptyset \Leftrightarrow d  \equiv 1 (\mod 8)$, \vspace{1.5ex}
\item[(2)] $C_{d-}(\mathbb{Q}_{p_i}) \neq \emptyset \Leftrightarrow (\frac{d}{p_i}) \ =\ 1 , \forall p_i \nmid d$,
\item[(3)] $C_{d-}(\mathbb{Q}_{p_i}) \neq \emptyset \Leftrightarrow (\frac{2D/(dp_i^{e^i-1})}{p_i}) \ =\ 1 , \forall  p_i \mid d$.
\end{enumerate}
}

\lemm{\label{lem2}
	Let $C_{d-}'$ be defined as above. Then\vspace{1.5ex}

\begin{enumerate}
\item[(1)] $C_{d-}'(\mathbb{Q}_2) \neq \emptyset \Leftrightarrow d  \equiv 1 ( \mod 8) \ \ or \ d - 2D/d\equiv 1 ( \mod 8)  $, \vspace{1ex}
\item[(2)] $C_{d-}'(\mathbb{Q}_{p_i}) \neq \emptyset \Leftrightarrow (\frac{d}{p_i}) \ =\ 1 , \forall p_i \nmid d$,
\item[(3)] $C_{d-}'(\mathbb{Q}_{p_i}) \neq \emptyset \Leftrightarrow (\frac{-2D/(dp_i^{e_i-1})}{p_i}) \ =\ 1 , \forall  p_i \mid d$.
\end{enumerate}
}

\lemm{\label{lem3}
Let $C_{d+}$ be defined as above. Then\vspace{1.5ex}

\begin{enumerate}
\item[(1)] $C_{d+}(\mathbb{Q}_2) \neq \emptyset \Leftrightarrow d  \equiv 1 (\mod 8)$, \vspace{1.5ex}
\item[(2)] $C_{d+}(\mathbb{Q}_{p_i}) \neq \emptyset \Leftrightarrow (\frac{d}{p_i}) \ =\ 1 , \forall p_i \nmid d$,
\item[(3)] $C_{d+}(\mathbb{Q}_{p_i}) \neq \emptyset \Leftrightarrow (\frac{-2D/(d p_i^{e^i-1})}{p_i}) \ =\ 1 , \forall  p_i \mid d$.
\end{enumerate}
}

\lemm{\label{lem4}
	Let $C_{d+}'$ be defined as above. Then\vspace{1.5ex}

\begin{enumerate}
\item[(1)] $C_{d+}'(\mathbb{Q}_2) \neq \emptyset \Leftrightarrow d \equiv 1( \mod 8) \ \ or \ d + 2D/d\equiv\ 1\ ( \mod\ 8)  $, \vspace{1ex}
\item[(2)] $C_{d+}'(\mathbb{Q}_{p_i}) \neq \emptyset \Leftrightarrow (\frac{d}{p_i}) \ =\ 1 , \forall p_i \nmid d$,
\item[(3)] $C_{d+}'(\mathbb{Q}_{p_i}) \neq \emptyset \Leftrightarrow (\frac{2D/(dp_i^{e_i-1})}{p_i}) \ =\ 1 , \forall  p_i \mid  d$.
\end{enumerate}
}
Since the Legendre symbol has the properties $(\frac{q^2}{p}) = 1$ and $p^2 \equiv 1 (\mod 8)$ for any $q \in \QQ$ and $p$ odd prime satisfying $\gcd (p,q) = 1$, there is no difference between $e_i = 1$ or $3$ in Lemmas \ref{lem1} to \ref{lem4}. Thus we assume $e_i = 1$, $\forall 1 \le i \le n$, or in other words, $D$ is squrefree, in the following part.

\section{Proof of the relation}\label{sec_pr}

From the Lemmas \ref{lem1} to \ref{lem4}, we see that whether an element is in $s_{\pm}(\phi)$ and $s'_{\pm}(\phi)$ or not can be determined by the congruence conditions and the Legendre symbol conditions. In this section, we give the proof of the main theorems.
In our proof, we show that the congruence conditions here can switch to some other Legendre symbol conditions. With the property of Legendre symbol we mentioned in section \ref{sec_2}, we will use $d$ instead of $\widetilde{d} = d \mod  \QQ^{\times 2}$ for convenience. Similarly, when we write $d_1d_2$ in the following statement, actually we talk about $\widetilde{d_1d_2} = d_1*d_2 \mod  \QQ^{\times 2}$.

For any $d \in \mathbb{Q}_{2D}$, let $g(d) = (g_{1}(d),g_{2}(d),\cdots,g_{n}(d))$. Similarly, for any $d' \in \mathbb{Q}_{-2D}$,
let $f(d') = (f_{1}(d'),f_{2}(d'),\cdots,f_{n}(d'))$, where
\begin{equation}\label{eq:array}
g_{i}(d) =\left\{
\begin{array}{ll}
(\frac{d}{p_i})        &\forall p_{i} \nmid d,\\
(\frac{2D/d}{p_i})      &\forall p_{i} \mid d.\\
\end{array}
\right.
\hspace{2em}
f_{i}(d') =\left\{\begin{array}{ll}
(\frac{d'}{p_i})        &\forall p_{i} \nmid d',\\
(\frac{-2D/d'}{p_i})      &\forall p_{i} \mid d'.\\
\end{array}\right.
\end{equation}
Particularly, we have $f_{i}(-1) = (\frac{-1}{p_i})=(-1)^{\frac{p_i - 1}{2}}$, $f_{i}(-D) = (\frac{2}{p_i}) = (-1)^{\frac{p_i^2 - 1}{8}}$, and $f_{i}(p_j) = (\frac{p_j}{p_i})$. It follows that 
\begin{equation} \label{3.2}
  f_{i}(p_j)f_{j}(p_i) =(\frac{p_j}{p_i})(\frac{p_i}{p_j}) = (-1)^{\frac{ p_i - 1}{2}\cdot \frac{ p_j - 1}{2}}.  
\end{equation}

Besides, we give the  following lemma to show that $g$ and $f$ are group homomorphisms:

\lemm{\label{lemma_1} For any $d_1,d_2 \in \QQ_{2D}$, and $d'_1,d'_2 \in \QQ_{-2D}$,
$$g(d_1)g(d_2) = g(d_1d_2),\  \forall d_1,d_2 \in \mathbb{Q}_{2D},$$
$$f(d'_1)f(d'_2) = f(d'_1d'_2),\ \forall d'_1,d'_2 \in \mathbb{Q}_{-2D}.$$
}
\pf{For any $d_1,d_2 \in \mathbb{Q}_{2D}$, we only need to verify it in $i$-th component. Now we divide the proof into four cases.

(1) For any $i$ such that $p_i \nmid d_1, d_2$, since $(\frac{d^2}{p_i}) = 1$,
it is obvious that
 $$g_{i}(d_1)g_{i}(d_2) =(\frac{d_1}{p_i})(\frac{d_2}{p_i}) = (\frac{d_1d_2}{p_i}) = g_{i}(d_1d_2).$$

(2) For any $i$ such that $p_i \mid d_1, d_2,$ we have $p_i \nmid {d_1d_2}$, then
$$g_{i}(d_1)g_{i}(d_2) =(\frac{2D/d_1}{p_i})(\frac{2D/d_2}{p_i}) = (\frac{4D^2/d_1d_2}{p_i}) =(\frac{d_1d_2 \mod \QQ^{\times 2}}{p_i}) = g_{i}(d_1d_2).$$

(3) For any $i$ such that $p_i \mid d_1, p_i \nmid d_2,$ we have $p_i \mid {d_1d_2}$.
From case(2) we get
$$g_{i}(d_1)g_{i}(d_1d_2) = g_{i}(d_1d_1d_2) =  g_{i}(d_2),$$
$$g_{i}(d_1d_2) =g^{2}_{i}(d_1)g_{i}(d_1d_2) = g_{i}(d_1)g_{i}(d_2).$$

(4) Similarly, we get the conclusion when  $p_i \nmid d_1, p_i \mid d_2$.

To summarize,
$$g_{i}(d_1)g_{i}(d_2) = g_{i}(d_1d_2).$$

Similarly, we have
$$f_{i}(d'_1)f_{i}(d'_2) = f_{i}(d'_1d'_2) .$$
}
In order to simplify the proof, we introduce $x(d) $ and $ y(d')$ here.

\defi{
Define $x(d) = (x_{1}(d),x_{2}(d),\cdots,x_{n}(d))$ and $y(d') = (y_{1}(d'),y_{2}(d'),\cdots,y_{n}(d')) \in \FF_2^n$,  such that $g_{i}(d) = (-1)^{x_{i}(d)}, f_{i}(d') = (-1)^{y_{i}(d')}, \ \forall d \in \QQ_{2D}, d'\in \QQ_{-2D}$. In other words,
\begin{equation}\label{eq:array2}
x_{i}(d) =\left\{
\begin{array}{ll}
0       &\forall g_{i}(d)  = 1,\\
1      &\forall g_{i}(d) = -1.\\
\end{array}
\right.
\hspace{2em}
y_{i}(d') =\left\{\begin{array}{ll}
0       &\forall f_{i}(d')  = 1,\\
1      &\forall f_{i}(d') = -1.\\
\end{array}\right.
\end{equation}
}

Then by Lemma \ref{lemma_1}, we have
$$x_{i}(d_1) + x_{i}(d_2) = x_{i}(d_1d_2),  \forall d_1,d_2 \in \mathbb{Q}_{2D},$$
and
$$y_{i}(d'_1)+ y_{i}(d'_2) = y_{i}(d'_1d'_2),  \forall d'_1,d'_2 \in \mathbb{Q}_{-2D}.$$
Thus $x$ and $y$ are group homomorphisms as well.

By the properties of Legendre symbols and the definition of $x_{i}(d_1)$ and $y_{i}(d'_1)$, the following proposition is straightforward. For example, the third one is straight forward corollary by \ref{3.2}.
\prop{\label{prop1} We have

(1) $x_{i}(d) = y_{i}(d),\quad \forall d \in \QQ_{2D}, 1 \leq i \leq n,$ such that $ p_i \nmid d$,

(2) $x_{i}(d) = y_{i}(d) + y_{i}(-1) =  y_{i}(-d),\quad \forall d \in \QQ_{2D}, 1 \leq i \leq n$, such that $p_i \mid d$,

(3) $x_{j}(p_i) + x_{i}(p_j) = y_{j}(p_i) + y_{i}(p_j) =y_{i}(-1)y_{j}(-1) ,\quad \forall 1 \leq i \neq j \leq n$.
}

\lemm{\label{lemma_2} Assume  $d = p_{i_1}p_{i_2} \cdots p_{i_k}$, then we have the following results:
$$ y_{i_1}(-1)+ \cdots + y_{i_k}(-1) = 0 \Leftrightarrow (\frac{-1}{p_{i_1}})\cdots (\frac{-1}{p_{i_k}}) = 1 \Leftrightarrow d \equiv 1,5 \mod 8,$$
$$y_{i_1}(-D)+ \cdots + y_{i_k}(-D) = 0 \Leftrightarrow (\frac{2}{p_{i_1}})\cdots (\frac{2}{p_{i_k}}) = 1 \Leftrightarrow d \equiv 1,7 \mod 8,$$
$$y_{i_1}(D)+ \cdots + y_{i_k}(D) = 0 \Leftrightarrow (\frac{-2}{p_{i_1}})\cdots (\frac{-2}{p_{i_k}})= 1 \Leftrightarrow d \equiv 1,3 \mod 8.$$
}

\pf{ Let $n_{i} = \# \{ p \mid p \text{ is a factor of } d,\  p \equiv i (\mod 8)\},\  i = 1,3,5,7$. Since $3^2 \equiv 5^2 \equiv 7^2 \equiv 1 \mod 8$, we can know the congruent number of d module 8 through the result for  $n_i$ module 2.
It's easy to know that
$$
\begin{aligned}
d \equiv 1 (\mod 8) & \Leftrightarrow n_3  \equiv n_5 \equiv n_7 (\mod 2), \\
d \equiv 3 (\mod 8) & \Leftrightarrow n_5 \equiv n_7 \equiv n_3+1(\mod 2), \\
d \equiv 5 (\mod 8) & \Leftrightarrow n_3 \equiv n_7 \equiv n_5+1 (\mod 2), \\
d \equiv 7 (\mod 8) & \Leftrightarrow n_3 \equiv n_5 \equiv n_7 +1(\mod 2). \\
\end{aligned}
$$
With the fact that
$$
\begin{aligned}
y_{i}(-1) = 0   &\Leftrightarrow \ (\frac{-1}{p_i}) = 1 &\Leftrightarrow p_i \equiv 1,5 (\mod 8), \\
y_{i}(-D) = 0   &\Leftrightarrow\ \ (\frac{2}{p_i}) = 1 &\Leftrightarrow p_i \equiv 1,7 (\mod 8), \\
y_{i}(D) = 0   &\Leftrightarrow\ (\frac{-2}{p_i}) = 1 &\Leftrightarrow p_i \equiv 1,3 (\mod 8),  \\
\end{aligned}
$$
we can get the result easily.
}

\thm{\label{thm1}Let $d = p_{i_1}p_{i_2} \cdots p_{i_k}$, $1 \le i_1 < i_2 \cdots < i_k \le n$. Then the condition $d \in s_-(\phi)$ is equivalent to $ \sum\limits_{j=1}^k y_{i_j}(p_m) = 0,\quad \forall 1 \le m \le n+1$.}
\pf{From Lemma \ref{lem1} we know
$$
\begin{aligned}
d = p_{i_1}p_{i_2} \cdots p_{i_k}  \in s_-(\phi) & \Leftrightarrow d \equiv 1 (\mod 8) \text{ and }  g_{m}(d) = 1, \forall 1\leq m \leq n, \\
&\Leftrightarrow d \equiv 1 (\mod 8) \text{ and } x_{m}(d) = 0, \forall 1\leq m \leq n .
\end{aligned}
$$
From Lemma \ref{lemma_2}, we have $d \equiv 1 \mod 8 \Leftrightarrow \delta(d) = \sum\limits_{j = 1}^k y_{i_j}(-1) = 0,$ and $\sum\limits_{j = 1}^k y_{i_j}(D) = 0$.
If we focus on the difference between $x_m(d)$ and the sum of $y_{i_j}(p_m)$, we have:
when $m \neq i_j, \ \forall 1 \le j \le k$:
	$$
	\begin{aligned}
	x_m(d) - \sum\limits_{j=1}^k y_{i_j}(p_m) & = \sum\limits_{j=1}^k (x_m(p_{i_j}) - y_{i_j}(p_m))\\
	& = \sum\limits_{j=1}^k (x_m(p_{i_j}) - x_{i_j}(p_m))\\
	& = \sum\limits_{j=1}^k (x_m(p_{i_j}) + x_{i_j}(p_m))\\
	& = \sum\limits_{j=1}^k (y_m(-1) y_{i_j}(-1))\\
	& = y_m(-1)\delta(d),
	\end{aligned}
	$$
when $ m = i_j, \ \exists 1 \le j \le k $:
	$$
	\begin{aligned}
	x_m(d) - \sum\limits_{j=1}^k y_{i_j}(p_m) & = \sum\limits_{j=1}^k (x_m(p_{i_j}) - y_{i_j}(p_m))\\
	& = \sum\limits_{1 \le j \le n ,m \neq i_j} (x_m(p_{i_j}) - x_{i_j}(p_m)) + y_{m}(-1)\\
	& = \sum\limits_{1 \le j \le n ,m \neq i_j} (x_m(p_{i_j}) + x_{i_j}(p_m))+ y_{m}(-1)\\
	& = \sum\limits_{1 \le j \le n ,m \neq i_j} (y_m(-1) y_{i_j}(-1)) + y_{m}(-1)\\
	& = y_m(-1)(\delta(d) - y_m(-1) + 1).
	\end{aligned}
	$$
It follows that
$$
\begin{aligned}
x_m(d) = \sum\limits_{j=1}^k y_{i_j}(p_m),\  \forall 1 \le m \le n & \Leftrightarrow \delta(d) = 0 \text{ or } y_m(-1) = 0,\  \forall 1 \le m \le n \\
& \Leftrightarrow \delta(d) = 0,
\end{aligned}
$$
and $\delta(d) = 0$ can be deduced from the conditions $d \equiv 1 \mod 8$ and $\sum\limits_{j=1}^k y_{i_j}(d') = 0, \ \forall d' \in \QQ_{-2D}$. Thus under the condition  $d \equiv 1 \mod 8$ or $\sum\limits_{j=1}^k y_{i_j}(p_{n+1}) = 0$, we have $x_m(d) = 0 \Leftrightarrow \sum\limits_{j=1}^k y_{i_j}(p_m) = 0,\  \forall 1 \le m \le n$.
And if $\sum\limits_{j=1}^k y_{i_j}(p_m) = 0,\  \forall 1 \le m \le n$, we have $$\sum\limits_{j = 1}^k y_{i_j}(D) = \sum\limits_{j = 1}^k \sum\limits_{m = 1}^n y_{i_j}(p_m) = \sum\limits_{m = 1}^n \sum\limits_{j = 1}^k y_{i_j}(p_m)= 0 .$$
Thus $d \in s_-(\phi)$ is equivalent to $ \sum\limits_{j=1}^k y_{i_j}(p_m) = 0, \forall 1 \le m \le n+1$.
}

\thm{\label{thm2}For any $d' \in \QQ_{-2D}$, $d' = p_{i_1}p_{i_2}\cdots p_{i_k} , \ 1\le i_1 < i_2 < \cdots < i_k \le n+1, \ k \ge 1$, we have $d' \in s_-'(\phi)$ if and only if $y_{j}(d') =  y_{j}(p_{i_1}) + y_{j}(p_{i_2}) + \cdots + y_{j}(p_{i_k}) = 0 \ \forall 1\le j \le n$.}

\pf{
	From Lemma \ref{lem2}, we have
$$
\begin{aligned}
d' = p_{i_1}p_{i_2} \cdots p_{i_k}  \in s'_-(\phi) & \Leftrightarrow d' \equiv 1 (\mod 8) \text{ or } d' - 2D/d' \equiv 1 (\mod 8)\text{ and }  f_{m}(d') = 1, \forall 1\leq m \leq n, \\
& \Leftrightarrow d' \equiv 1 (\mod 8) \text{ or } d' - 2D/d' \equiv 1 (\mod 8)\text{ and } y_{m}(d') = 0, \forall 1\leq m \leq n .
\end{aligned}
$$
Thus the forward direction is obvious, let us turn to the backward direction.

(\romannumeral1) If $i_k \le n$, without loss of generality, we may assume $d' = p_{1}p_{2}\cdots p_{k}$. We have
	$$ \sum\limits_{i = 1}^{k} y_{j}(p_i) = 0,\ \forall  k < j \le n,$$
	and
	$$  y_{i}(D) +\sum\limits_{j = k+1}^{n} y_{i}(p_j) = 0, \   \forall 1 \le i \le k. $$
	Then we get
	$$
	\begin{aligned}
	0 &= \sum\limits_{i = 1}^{k} (y_{i}(D) +\sum\limits_{j = k+1}^{n} y_{i}(p_j) ) + \sum\limits_{j = k+1}^{n}\sum\limits_{i = 1}^{k}y_{j}(p_i)  \\
	 &= \sum\limits_{j = k+1}^{n}\sum\limits_{i = 1}^{k}(y_{j}(p_i) + y_{i}(p_j)) + \sum\limits_{i = 1}^{k} y_{i}(D)\\
	 &= \sum\limits_{j = k+1}^{n}\sum\limits_{i = 1}^{k} y_{i}(-1)y_{j}(-1) + \sum\limits_{i = 1}^{k} y_{i}(D)\\
	 &= (\sum\limits_{j = k+1}^{n} y_{j}(-1))(\sum\limits_{i = 1}^{k}y_{i}(-1)) + \sum\limits_{i = 1}^{k} y_{i}(D).
	\end{aligned}
	$$

	 Let $e = p_{k+1}p_{k+2}\cdots p_{n}$ , $de = D$. Then we distinguish among three cases.

	Case 1: $\sum\limits_{i = 1}^{k}y_{i}(-1) = \sum\limits_{i = 1}^{k} y_{i}(D) = 0 $. By lemma\ref{lemma_2}, in this case we have $d'$ satisfies $d' \equiv 1,5 (\mod 8)$ and $d' \equiv 1,3 (\mod 8)$, so $d' \equiv 1 (\mod 8) $.

	Case 2: $\sum\limits_{i = 1}^{k}y_{i}(-1) = 1, \sum\limits_{j = k+1}^{n} y_{j}(-1)  = \sum\limits_{i = 1}^{k} y_{i}(D) = 0 $. By lemma\ref{lemma_2}, in this case we have $d'$ satisfies $d' \equiv 3,7 (\mod 8)$ and $d' \equiv 1,3 (\mod 8)$, $e$ satisfies  $e \equiv 1,5 (\mod 8)$. Hence $d' - 2D/d' \equiv d - 2e \equiv 1 (\mod 8)$.

	Case 3: $\sum\limits_{j = k+1}^{n} y_{j}(-1) = \sum\limits_{i = 1}^{k}y_{i}(-1) = \sum\limits_{i = 1}^{k} y_{i}(D) = 1$.  By lemma\ref{lemma_2}, in this case we have $d'$ satisfies  $d' \equiv 3,7 (\mod 8)$ and $d' \equiv 5,7 (\mod 8)$, $e$ satisfies  $e \equiv 3,7 (\mod 8)$. Hence $d' - 2D/d' \equiv d' - 2e \equiv 1 (\mod 8)$.

	(\romannumeral2) If $i_k = n+1$, without loss of generality, assume $d' = -p_{1}p_{2}\cdots p_{k-1}$. We have
	$$ y_{j}(-1) + \sum\limits_{i = 1}^{k-1} y_{j}(p_i) = 0, \ \forall  k \le j \le n,$$
	and
	$$  y_{i}(-D) +\sum\limits_{j = k}^{n} y_{i}(p_j) = 0 , \   \forall 1 \le i \le k-1. $$
Then we get
	$$
	\begin{aligned}
	0 &= \sum\limits_{i = 1}^{k-1} (y_{i}(-D) +\sum\limits_{j = k}^{n} y_{i}(p_j) ) + \sum\limits_{j = k}^{n}(\sum\limits_{i = 1}^{k} y_{j}(p_i) + y_{j}(-1))  \\
	 &= \sum\limits_{j = k}^{n}\sum\limits_{i = 1}^{k-1}(y_{j}(p_i) + y_{i}(p_j)) + \sum\limits_{i = 1}^{k-1} y_{i}(-D) + \sum\limits_{j = k}^{n} y_{j}(-1)\\
	 &= \sum\limits_{j = k}^{n}\sum\limits_{i = 1}^{k-1} y_{i}(-1)y_{j}(-1) + \sum\limits_{i = 1}^{k-1} y_{i}(-D)+ \sum\limits_{j = k}^{n} y_{j}(-1)\\
	&= (\sum\limits_{j = k}^{n} y_{j}(-1))(1+\sum\limits_{i = 1}^{k-1}y_{i}(-1)) + \sum\limits_{i = 1}^{k-1} y_{i}(-D).
	\end{aligned}
	$$
	Let $e = p_{k}p_{k+1}\cdots p_{n}$ , $de = -D$. Also, we distinguish among three cases.

	Case 1: $\sum\limits_{i = 1}^{k-1}y_{i}(-1) = 1$ and  $\sum\limits_{i = 1}^{k} y_{i}(-D) = 0 $. By lemma\ref{lemma_2}, in this case we have $-d' \equiv 3,7 (\mod 8)$ and $d' \equiv 1,7 (\mod 8)$, so $-d' \equiv 7 (\mod 8) $ and $d' \equiv 1 (\mod 8) $

	Case 2: $\sum\limits_{i = 1}^{k-1}y_{i}(-1) = 0, \sum\limits_{j = k}^{n} y_{j}(-1)  = \sum\limits_{i = 1}^{k-1} y_{i}(-D) = 0 $. By lemma\ref{lemma_2}, in this case we have $-d' \equiv 1,5 (\mod 8)$ and $-d' \equiv 1,7 (\mod 8)$, $e$ satisfies  $e \equiv 1,5 (\mod 8)$. Hence $d' - 2D/d' \equiv d + 2e \equiv 1 (\mod 8)$.

	Case 3: $\sum\limits_{i = 1}^{k-1}y_{i}(-1) = 0$, $\sum\limits_{j = k}^{n} y_{j}(-1) =  \sum\limits_{i = 1}^{k} y_{i}(D) = 1$. By lemma\ref{lemma_2}, in this case we have $-d' \equiv 1,5 (\mod 8)$ and $-d' \equiv 3,5 (\mod 8)$, $e$ satisfies  $e \equiv 3,7 (\mod 8)$. Hence $d' - 2D/d' \equiv d + 2e \equiv 1 (\mod 8)$.

	To sum up, $d' \in s_-'(\phi)$ is equivalent to $y_{j}(d') = 0 \ \forall 1\le j \le n$.
}
\thm{\label{thm3}For a fixed $d_1 = p_{i_1}p_{i_2} \cdots p_{i_k}, 1\le i_1 < \cdots < i_k \le n$, then one and only one $d = \sigma d_1 \in s_+(\phi)$ is equivalent to $ \sum\limits_{j=1}^k x_{i_j}(p_m) = 0, \forall 1 \le m \le n$.
}
\pf{From Lemma \ref{lem1} we know
$$
\begin{aligned}
d = \sigma p_{i_1}p_{i_2} \cdots p_{i_k}  \in s_+(\phi) & \Leftrightarrow d \equiv 1 (\mod 8) \ and \  f_{m}(d) = 1, \forall 1\leq m \leq n, \\
&\Leftrightarrow d \equiv 1 (\mod 8) \ and  \ y_{m}(d) = 0, \forall 1\leq m \leq n .
\end{aligned}
$$
By Lemma \ref{lemma_2}, we have $d_1 \equiv \pm 1 \mod 8 \Leftrightarrow \sum\limits_{j = 1}^k x_{i_j}(D) = 0$.
If we focus on the difference between $y_m(d)$ and the sum of $x_{i_j}(p_m)$, we have:
when $m \neq i_j, \ \forall 1 \le j \le k$:
	$$
	\begin{aligned}
	y_m(d) - \sum\limits_{j=1}^k x_{i_j}(p_m) & = y_m(\sigma) + \sum\limits_{j=1}^k (y_m(p_{i_j}) - x_{i_j}(p_m))\\
	& =  y_m(\sigma) +\sum\limits_{j=1}^k (x_m(p_{i_j}) - x_{i_j}(p_m))\\
	& =  y_m(\sigma) +\sum\limits_{j=1}^k (x_m(p_{i_j}) + x_{i_j}(p_m))\\
	& =  y_m(\sigma) +\sum\limits_{j=1}^k (y_m(-1) y_{i_j}(-1))\\
	& =  y_m(\sigma) +y_m(-1)\delta(d_1),
	\end{aligned}
	$$
when $ m = i_j, \ \exists 1 \le j \le k $:
	$$
	\begin{aligned}
	y_m(d) - \sum\limits_{j=1}^k x_{i_j}(p_m) & =  y_m(\sigma) + \sum\limits_{j=1}^k (y_m(p_{i_j}) - x_{i_j}(p_m))\\
	& =  y_m(\sigma) + \sum\limits_{1 \le j \le n ,m \neq i_j} (x_m(p_{i_j}) - x_{i_j}(p_m)) + y_{m}(-1)\\
	& = y_m(\sigma) +  \sum\limits_{1 \le j \le n ,m \neq i_j} (x_m(p_{i_j}) + x_{i_j}(p_m))+ y_{m}(-1)\\
	& =  y_m(\sigma) + \sum\limits_{1 \le j \le n ,m \neq i_j} (y_m(-1) y_{i_j}(-1)) + y_{m}(-1)\\
	& =  y_m(\sigma) + y_m(-1)(\delta(d_1) - y_m(-1) + 1).
	\end{aligned}
	$$

First we prove the forward direction.
In fact if $\sigma = 1$, then by $d = \sigma d_1 \equiv 1 \mod 8$, we have $y_m(\sigma) = \delta(d_1) = 0$. As a result, $0 = y_m(d) = \sum\limits_{j=1}^k x_{i_j}(p_m),\  \forall 1\le m \le n$. Otherwise, we know $\sigma = -1$,
similarly, we have $y_m(\sigma) =  y_m(-1)$ and $\delta(d_1) = 1$. We can also get $0 = y_m(d) = \sum\limits_{j=1}^k x_{i_j}(p_m),\  \forall 1\le m \le n$. Our result follows.

For the backward direction, if $\sum\limits_{j=1}^k x_{i_j}(p_m) = 0, \forall 1 \le m \le n$, we have:
$$0 = \sum\limits_{m=1}^n \sum\limits_{j=1}^k x_{i_j}(p_m) = \sum\limits_{j=1}^k\sum\limits_{m=1}^n x_{i_j}(p_m) =\sum\limits_{j=1}^k x_{i_j}(D) = \sum\limits_{j=1}^k y_{i_j}(-D). $$
By Lemma \ref{lemma_2}, we have $d_1 \equiv 1,7 (\mod 8)$. Select $d = d_1$ when $d_1 \equiv 1 (\mod 8)$, and $d = -d_1$ when $d_1 \equiv 7 (\mod 8)$. Then $d \equiv 1 (\mod 8)$. Also we have $y_m(\sigma) = 0,\  \delta(d_1) = 0$ or $y_m(\sigma) =  y_m(-1),\ \delta(d_1) = 1$. Thus $0 = \sum\limits_{j=1}^k x_{i_j}(p_m)= y_m(d),\  \forall 1\le m \le n$. Since $d \equiv 1 (\mod 8)$ and $y_m(d) = 0, \forall 1\le m \le n$, we know $d \in s_+(\phi)$.
}

\thm{\label{thm4}For any $d \in \QQ_{2D}$, $d = p_{i_1}p_{i_2}\cdots p_{i_k} , \ 1\le i_1 < i_2 < \cdots < i_k \le n, \ k \ge 1$, we have $d \in s_+'(\phi)$ if and only if $x_{j}(d) =  x_{j}(p_{i_1}) + x_{j}(p_{i_2}) + \cdots + x_{j}(p_{i_k}) = 0 ,\ \forall 1\le j \le n$.}

\pf{
	From Lemma\ref{lem2}, we have
$$
\begin{aligned}
d = p_{i_1}p_{i_2} \cdots p_{i_k}  \in s_+'(\phi) & \Leftrightarrow d \equiv 1 (\mod 8) \text{ or } d + 2D/d \equiv 1 (\mod 8)\text{, and }  g_{m}(d) = 1, \forall 1\leq m \leq n, \\
& \Leftrightarrow d \equiv 1 (\mod 8) \text{ or } d + 2D/d \equiv 1 (\mod 8)\text{, and } x_{m}(d) = 0, \forall 1\leq m \leq n .
\end{aligned}
$$
Thus the forward direction is obvious, let us turn to the backward direction. Without loss of generality, assume $d = p_{1}p_{2}\cdots p_{k}$. We have
	$$ \sum\limits_{i = 1}^{k} x_{j}(p_i) = 0,\ \forall  k < j \le n,$$
	and
	$$  x_{i}(D) +\sum\limits_{j = k+1}^{n} x_{i}(p_j) = 0, \   \forall 1 \le i \le k. $$
	Then we get
	$$
	\begin{aligned}
	0 &= \sum\limits_{i = 1}^{k} (x_{i}(D) +\sum\limits_{j = k+1}^{n} x_{i}(p_j) ) + \sum\limits_{j = k+1}^{n}\sum\limits_{i = 1}^{k}x_{j}(p_i)  \\
	 &= \sum\limits_{j = k+1}^{n}\sum\limits_{i = 1}^{k}(y_{j}(p_i) + y_{i}(p_j)) + \sum\limits_{i = 1}^{k} x_{i}(D)\\
	 &= \sum\limits_{j = k+1}^{n}\sum\limits_{i = 1}^{k} y_{i}(-1)y_{j}(-1) + \sum\limits_{i = 1}^{k} y_{i}(-D)\\
	 &= (\sum\limits_{j = k+1}^{n} y_{j}(-1))(\sum\limits_{i = 1}^{k}y_{i}(-1)) + \sum\limits_{i = 1}^{k} y_{i}(-D)
	\end{aligned}
	$$

	 Let $e = p_{k+1}p_{k+2}\cdots p_{n}$ , $de = D$. There are three cases to check.

	Case 1: $\sum\limits_{i = 1}^{k}y_{i}(-1) = \sum\limits_{i = 1}^{k} y_{i}(-D) = 0 $. By Lemma \ref{lemma_2}, in this case we have $d$ satisfies $d \equiv 1,5 (\mod 8)$ and $d \equiv 1,7 (\mod 8)$, so $d \equiv 1 (\mod 8) $.

	Case 2: $\sum\limits_{i = 1}^{k}y_{i}(-1) = 1, \sum\limits_{j = k+1}^{n} y_{j}(-1)  = \sum\limits_{i = 1}^{k} y_{i}(-D) = 0 $. By Lemma \ref{lemma_2}, in this case we have $d$ satisfies $d \equiv 3,7 (\mod 8)$ and $d \equiv 1,7 (\mod 8)$, $e$ satisfies  $e \equiv 1,5 (\mod 8)$. Hence $d + 2D/d \equiv d + 2e \equiv 1 (\mod 8)$.

	Case 3: $\sum\limits_{j = k+1}^{n} y_{j}(-1) = \sum\limits_{i = 1}^{k}y_{i}(-1) = \sum\limits_{i = 1}^{k} y_{i}(-D) = 1$.  By Lemma \ref{lemma_2}, in this case we have $d$ satisfies  $d \equiv 3,7 (\mod 8)$ and $d \equiv 3,5 (\mod 8)$, $e$ satisfies  $e \equiv 3,7 (\mod 8)$. Hence $d + 2D/d \equiv d + 2e \equiv 1 (\mod 8)$.

	To sum up, $d \in s_+'(\phi)$ is equivalent to $x_{j}(d) = 0 \ \forall 1\le j \le n$.
}

\defi{\label{def3}Let $l(1) = 0$, and $l(d) = \max\{1 \le i \le n\mid p_i \text{ divides } d\}$, for any $d \in \QQ_{-2D}$, $d \neq 1$. Similarly, let $l'(1) = 0$, $p_{n+1} = -1, $ and $l'(d) = i_k,$ for $d' = p_{i_1}p_{i_2}\cdots p_{i_k} \in \QQ_{-2D}$, $1 \le i_1 < i_2 < \cdots < i_k \le n+1 ,$ $ k \geq 1$. }

In the following, we give the proof of Theorem \ref{thm_1}.
\thm{ \label{thm_n1} Let $E_-$, $E_-'$ be the elliptic curves defined above. Then     $$2^{-t_{2D}} = T_{2D} = \# S^{(\phi)}(E_-)/\# S^{(\phi)}(E'_{-}) = 1/2.$$ 
}
\pf{
 Let us define a matrix $Y = [m'_{i,j}]_{(n+1)\times n}$, where $m'_{i,j} = y_{j}(p_i)$.
 Since $y$ is a homomorphism, $Y$ encodes the information of $\{y_{i}(d)\}$. Let $r'_{i},1 \le i \le n+1$ be the row vectors of $Y$ and $c'_{j},1\le j \le n$ the column ones as below.
\begin{equation}\label{eq:array3}
\bordermatrix{%
       & c'_1       & c'_2     &\cdots     &c'_n\cr
r'_1    &y_{1}(p_1) & y_{2}(p_1) & \cdots  &y_{n}(p_1)\cr
\vdots & \vdots    &\vdots   &\cdots     &\vdots\cr
r'_n 		& y_{1}(p_n) & y_{2}(p_n) & \cdots  &y_{n}(p_n)\cr
r'_{n+1} 		&y_{1}(p_{n+1}) & y_{2}(p_{n+1}) & \cdots  &y_{n}(p_{n+1})
}
\begin{array}{c}
      \longleftarrow  \text{corresponding to $d' = p_1$} \\
		\vdots  \\
	  \longleftarrow  \text{corresponding to $d' = p_n$}\\
	  \longleftarrow  \text{corresponding to $d' = -1$}\\
    \end{array}
\end{equation}

Since $y$ is a homomorphism, to fulfill the condition $y_{i_1}(d')+y_{i_2}(d')+\cdots +y_{i_k}(d') = 0, \ \forall d' \in \QQ_{-2D}$ in Theorem \ref{thm1}, all we need is $y_{i_1}(d')+y_{i_2}(d')+\cdots +y_{i_k}(d') = 0, \quad \forall d' = p_{i}, \quad 1 \le i \le n+1,$ which is equivalent to the linear dependence condition $c'_{i_1}+c'_{i_2}+\cdots+c'_{i_k} = 0$.

Similarly, for a fixed $d' = p_{i_1}p_{i_2}\cdots p_{i_k} \neq 1, 1\le i_1 < i_2 < \cdots < i_k \le n+1$,  from Theorem \ref{thm2}, we have
$$
	\begin{aligned}
	d' \in s'(\phi) &\Leftrightarrow y(d') = 0 \\
					 &\Leftrightarrow r'_{i_1}+r'_{i_2}+\cdots +r'_{i_k}=0 \quad (y \text{ is a homomorphism}).
	\end{aligned}
$$
 Let $l(d), l'(d')$ be as in the Definition \ref{def3}, $R' = \{i>0 \mid \exists d \in s_-(\phi) , i = l(d) \}$ and $C' = \{i>0 \mid \exists d' \in s'_-(\phi) , i = l'(d') \}$.
 As the rank of $Y$ is a certain number for a given $Y$, we have $\# C' = \# R' + 1 = n+1-\rank(Y)$. Besides, since for any $ 0 \le i \le n$, either
 $$\#\{ d \mid d \in s_-(\phi) ,  l(d) \le i + 1\} = \#\{ d \mid d \in s_-(\phi) , l(d) \le i \}  $$ or
$$\#\{ d \mid d \in s_-(\phi) , l(d) \le i+1 \} = 2 \#\{ d \mid d \in s_-(\phi) , l(d) \le i\} , $$
we have $\# R' = \log_{2}(\# s_-(\phi))$, and similarly, $\# C' = \log_{2}(\# s'_-(\phi))$. In this way, we have  $\# s'_-(\phi) = 2\# s_-(\phi)$ immediately. By the fact $\#  S^{(\phi)}(E_-) = 2 \# s_-(\phi)$ and $\# S^{({\phi})}({E_-'}) = 2 \# s_-'(\phi)$, we have $\# S^{({\phi})}({E_-'}) = 2 \#  S^{(\phi)}(E_-)$. It follows that 
$$2^{-t_{2D}} = T_{2D} = \# S^{(\phi)}(E_-)/\# S^{(\phi)}(E'_{-}) = 1/2.$$
}

Similarly, we have the proof of Theorem \ref{thm_2}.

\thm{ \label{thm_n2} Let $E_+$, $E_+'$ be the elliptic curves defined above. Then     $$2^{-t_{-2D}} = T_{-2D} = \# S^{(\phi)}(E_+)/\# S^{(\phi)}(E'_{+}) = 1.$$
}
\pf{
 Let us define Matrix $X = [m_{i,j}]_{n\times n}$, where $m_{i,j} = x_{j}(p_i)$, which encodes information of $\{x_{i}(d)\}$. Let $r_{i},1 \le i \le n+1$ be the row vectors of $X$ and $c_{j},1\le j \le n$ be the column vectors.
\begin{equation}\label{eq:array4}
\bordermatrix{%
       & c_1       & c_2     &\cdots     &c_n\cr
r_1    &x_{1}(p_1) & x_{2}(p_1) & \cdots  &x_{n}(p_1)\cr
r_2    &x_{1}(p_2) & x_{2}(p_2) & \cdots  &x_{n}(p_2)\cr
\vdots & \vdots    &\vdots   &\cdots     &\vdots\cr
r_n 		& x_{1}(p_n) & x_{2}(p_n) & \cdots  &x_{n}(p_n)
}
\begin{array}{c}
      \longleftarrow  \text{corresponding to $d = p_1$} \\
     \longleftarrow  \text{corresponding to $d = p_2$} \\
		\vdots  \\
	  \longleftarrow  \text{corresponding to $d = p_n$}\\
    \end{array}
\end{equation}
Similarly, we have the condition $d = p_{i_1}p_{i_2} \cdots p_{i_k}  \in s_+(\phi)$ is equivalent to $r_{i_1}+r_{i_2}+\cdots +r_{i_k}=0$, And one of $ d $ and $-d \in s'_+(\phi)$ is equivalent to $c_{i_1}+c_{i_2}+\cdots +c_{i_k}=0$.

 Let $l(d)$ be as in the Definition \ref{def3}. Define $R= \{i>0\mid \exists d \in s_+(\phi) , i = l(d) \}$ and $C = \{i>0 \mid \exists d' \in s'_+(\phi) , i = l(d') \}$. Same as the Proof of Theorem \ref{thm_1}, we have
$$\# R = \# C = n-\rank(X), \quad \# R = \log_{2}(\# s_+(\phi)) \quad \# C = \log_{2}(\# s'_+(\phi)).$$ Similarly, we have  $\# s'_+(\phi) = \# s_+(\phi)$ immediately. By the fact $\#  S^{(\phi)}(E_+) = 2 \# s_+(\phi)$ and $\# S^{({\phi})}({E_+'}) = 2 \# s_+'(\phi)$, we have $\# S^{({\phi})}({E_+'}) = \#  S^{(\phi)}(E_+)$. Thus
$$2^{-t_{2D}} = T_{2D} = \# S^{(\phi)}(E_-)/\# S^{(\phi)}(E'_{-}) = 1.$$
}

\section{Application}

In this section, we give a method to quickly get the Mordell-Weil rank of part of Elliptic Curves in these forms.

In traditional method, determine the structure of $\phi$-Selmer group needs to traverse $\QQ(S,2)$, Which follows exponential computational complexity.
In the proof of our main theorem, the problem calculating $s_\pm(\phi)$ and $s'_\pm(\phi)$ are translated into the calculation of linear dependence of rows and columns in matrix $X$, $Y$ in $\FF_2$. Moreover, if we want to calculate the sizes of $s_\pm(\phi)$ and $s'_\pm(\phi)$, we only need to compute the ranks of $X$ and $Y$. Since by Gaussian elimination, we can determine the rank and the linear dependence relation of matrix with computational complexity $O(n^3)$ at most. We can quickly figure out the structure of $\phi$-Selmer groups of $E$.

In particular, if $\rank(X)= n$ (resp. $\rank(Y)= n$), we have $\# s'_+(\phi) = \# s_+(\phi) = 1$ (resp. $\# s'_-(\phi) =2 \# s_-(\phi) = 2$). It follows that $S^{(\phi)}(E_+)\cong S^{({\phi})}({E'_+}) \cong \ZZ/2\ZZ$ (and resp. $S^{({\phi})}({E'_-}) \cong (\ZZ/2\ZZ)^2 $, $s^{(\phi)}(E_+)\cong \ZZ/2\ZZ$.)
In addition, from, for example, \cite{gtm106}, we get
$$
\begin{aligned}
r_{E_{2D}} + \dim(\Sha(E_-)[{\phi}]) + \dim(\Sha(E'_-)[{\phi}]) &= \dim(S^{(\phi)}(E_-)) + \dim(S^{({\phi})}({E'_-})) - 2, \\
&= \dim(s_-(\phi)) + \dim(s'_-(\phi)),\\
&= \# R' + \# C',\\
&= 2n + 1 - 2\rank(Y).
\end{aligned}
$$
Similarly,
$$
\begin{aligned}
r_{E_{-2D}} + \dim(\Sha(E_+)[{\phi}]) + \dim(\Sha(E'_+)[{\phi}]) &= \dim(S^{(\phi)}(E_+)) + \dim(S^{({\phi})}({E'_+})) - 2, \\
&= \dim(s_+(\phi)) + \dim(s'_+(\phi)),\\
&= \# R + \# C,\\
&= 2n  - 2\rank(X).
\end{aligned}
$$
The $\dim$ here means $\dim_2$. Therefore, if $\rank(X) = n$, we have $r_{E_{+2D}} = 0$.
Besides, by Birch–Stephens \cite{parity}, the parity of $t_A$ is the same as that of the root number of $E_A$. 
$$(-1)^{t_A} = w(E_A).$$ 
Let $L_E(s)$ be the $L$ function associated to E, an elliptic curve over $\QQ$ of conductor $N$ and
$$
\Lambda_{E}(s) = (2\pi)^{-s}N(E)^{\frac{s}{2}}L_{E}(s).
$$
By Modularity theorem [3], we have $\Lambda_{E}(s) = W(E)\Lambda_{E}(2-s)$,
Then by $W(E_{-2D}) = -1$, we have $L_E(1) = 0$. By Birch and Swinnerton-Dyer conjecture \cite{gtm106}, we have $r_{E_{2D}} \geq 1$. Therefore, if $\rank(Y) = n$, we have $r_{E_{2D}} = 1$ and $\Sha_{E}[\phi] = \Sha_{E'}[{\phi}] = \emptyset$. Hence $E_{-2D}(\QQ) \cong E_{tors}\oplus\ZZ \cong \ZZ/2\ZZ \oplus \ZZ$.
Overall, we can quickly sieve some of the elliptic curves in this form which has Mordell-Weil rank $0$ or $1$.







	\bibliography{ref1}
	\bibliographystyle{amsplain}

\end{document}